\magnification=\magstep1
\input epsf

~~\vskip 1in

\centerline{\bf Surfaces, circles, and solenoids}

\vskip .2in

\centerline {\bf R. C. Penner}

\centerline{Departments of Mathematics and Physics/Astronomy}

\centerline{University of Southern California}

\centerline {Los Angeles, CA 90089}

\vskip .2in

\centerline{August 23, 2004}

\vskip .2in

\leftskip .5in\rightskip .5in

\noindent {\bf Abstract}~~Previous work of the author has developed coordinates on bundles over the classical Teichm\"uller spaces of punctured
surfaces and on the space of cosets of the M\"obius group in the group of orientation-preserving homeomorphisms of the circle, and this work is
surveyed here.  Joint work with Dragomir \v Sari\' c is also sketched which extends these results to the setting of the Teichm\"uller space of the
solenoid of a punctured surface, which is defined here in analogy to Dennis Sullivan's original definition for the case of closed surfaces.  Because
of relations with the classical modular group, the punctured solenoid and its Teichm\"uller theory have connections with
number theory.
 
\leftskip=0ex\rightskip=0ex

~~\vskip .2in

\centerline{ \bf Introduction}\vskip .2in

\noindent The ``lambda length'' of a pair of horocycles in upper halfspace centered at
$u,v\in{\bf R}$ of respective diameters
$c,d$ is given by
$\sqrt{2\over{cd}}~|u-v|$ and is roughly the hyperbolic distance between them.
(See $\S$1 for more precision.)  These invariants can be used to devise coordinates in several different guises: for the Teichm\"uller space of
punctured surfaces [7]; for the space of cosets of the M\"obius group of real fractional linear transformations in the
topological group of all orientation-preserving homeomorphisms of the circle, which forms a generalized universal Teichm\"uller space [8]; and for
the Teichm\"uller space of the ``punctured solenoid'', which is the punctured analogue introduced in [10] (and defined in $\S$4) of the space studied
by Sullivan [12] to analyze dynamical properties of the mapping class group actions on the Teichm\"uller spaces for closed surfaces. In fact in each
case, lambda lengths give coordinates for the ``decorated'' Teichm\"uller space rather than the Teichm\"uller space. (The respective notions of
decoration are defined in
$\S\S$ 2,3,4.)  Furthermore, the manifestation of lambda lengths as coordinates on the decorated Teichm\"uller space of the
punctured solenoid is the first step of a larger ongoing program with \v Sari\' c [10] to extend the decorated Teichm\"uller theory [7]-[9] to the
solenoid. 

\vskip .1in

\noindent To define the punctured solenoid ${\cal H}$ as a topological space, for definiteness fix the ``modular''  group $G=PSL_2({\bf Z})$ of
integral fractional linear transformations, let $\hat G$ denote its pro-finite completion (whose definition is recalled in $\S$4), let ${\bf D}$
denote the open unit disk in the complex plane, and define ${\cal H}=({\bf D}\times \hat{G})/G$, where $\gamma\in G$ acts on
$(z,t)\in{\bf D}\times\hat{G}$ by $\gamma (z,t)=(\gamma z,t\gamma ^{-1})$.  In analogy to the
case of punctured surfaces, we may produce appropriate geometric structures on ${\cal H}$ by taking suitable quotients $({\bf
D}\times\hat{G})/G$ by other actions of $G$ on ${\bf D}\times\hat{G}$.  As a pro-finite completion, the punctured solenoid itself is defined
essentially number theoretically in terms of finite-index subgroups of the modular group, and aspects of its Teichm\"uller theory 
bear close relation to classical questions in number theory (as mentioned at the end of $\S$5, which also contains other concluding and speculative
remarks).

\vskip .1in

\noindent We take this opportunity to correct Theorem~6.4 from [8].  See the remarks following Theorem~8 for the correction to the
universal Teichm\"uller theory and Proposition~12 for the corresponding affirmative statement for the solenoid.

\vskip .2in

\noindent {\bf Acknowledgement}~It was discussions with Mahmoud Zeinalian that led to the original idea of employing lambda lengths for solenoids,
and is  a pleasure to thank him, as well as Bob Guralnick for useful comments on classical number theory.  The new
material in
$\S$4 on the punctured solenoid is joint work [10] with Dragomir \v Sari\' c, who patiently explained his earlier work [11], and it
is also a pleasure to thank him for many stimulating conversations.

\vskip .3in

\noindent {\bf 1. Background}\vskip .2in

\noindent Define the Minkowski inner product $<\cdot ,\cdot >$ on ${\bf R}^3$ whose quadratic form
is given by
$x^2+y^2-z^2$ in the usual coordinates.   The upper sheet
$${\bf H}=\{ u=(x,y,z)\in{\bf R}^3:<u,u>=-1~{\rm and}~z>0\}$$ of the two-sheeted hyperboloid is isometric to the hyperbolic plane.
Indeed, identifying the Poincar\'e disk ${\bf D}$ with the open unit disk at height zero about the origin in ${\bf
R}^3$, central projection ${\bf H}\to {\bf D}$ from $(0,0,-1)\in{\bf R}^3$ establishes an isometry. 
Moreover, the open positive light cone
$$L^+=\{ u=(x,y,z)\in{\bf R}^3:<u,u> =0~{\rm and}~z>0\}$$ is identified with the collection of all horocycles in ${\bf H}$ via 
the affine duality $u\mapsto h(u)=\{ w\in{\bf H}:<w,u>=-1\}$.
Identifying the unit circle $S^1$ with the boundary of ${\bf D}$, the central projection extends continuously to the projection $\Pi:L^+\to S^1$
which maps a horocycle in
$L^+$ to its center in $S^1$.

\vskip .1in

\noindent 
Define a ``decorated geodesic'' to be an unordered pair $\{ h_0 ,h_1\}$ of horocycles with distinct centers
in the
hyperbolic plane, so there is a well-defined geodesic connecting the centers of $h_0$ and $h_1$; the two
horocycles may or may not be disjoint, and there is a well-defined signed hyperbolic distance $\delta$  between them
(taken to be positive if and only if 
$h_0\cap h_1=\emptyset$) as illustrated in the two cases of Figure~1.  The {\it lambda length} of the decorated geodesic $\{ h_0,h_1\}$ 
is defined to be the transform $\lambda (h_0,h_1)=\sqrt{2~{\rm exp}~\delta}$.  Taking this particular transform renders the
identification
$h$ geometrically natural in the sense that 
$\lambda (h(u_0),h(u_1))=\sqrt{-<u_0,u_1>}$, for $u_0,u_1\in L^+$ as one can check.

\vskip .2in

~~~{{{\epsffile{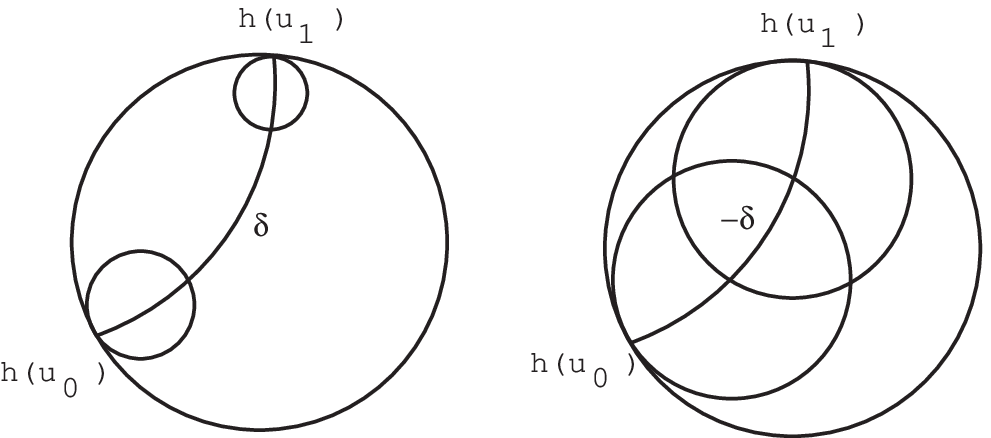}}}}

\vskip .1in

\centerline{{\bf Figure 1}~~{Decorated geodesics.}}

\vskip .2in

\noindent Three useful lemmas (with computational proofs, which we do not reproduce here) are as follows:

\vskip .1in 

\noindent {\bf Lemma 1} [7;Lemma 2.4] {\it Given three rays $\vec r_0,\vec r_1,\vec r_2\subseteq L^+$ from the origin which contain linearly
independent vectors and given three numbers $\lambda _0,\lambda _1,\lambda _2\in{\bf R}_+$, there are unique points $u_i\in \vec r_i$, for $i=0,1,2$,
so that
$\lambda (h(u_i),h(u_j))=\lambda _k$, where $\{ i,j,k\} =\{ 0,1,2\}$.  The points $u_0,u_1,u_2$ depend continuously on 
$\lambda _0,\lambda _1,\lambda
_2$ and on $\vec r_0,\vec r_1,\vec r_2$.} 

\vskip .1in  

\noindent {\bf Lemma 2} [7;Lemma 2.3] {\it Given two points $u_0,u_1\in L^+$, which do not lie on a common ray through the origin, and
given two numbers $\lambda _0,\lambda _1\in{\bf R}_+$, there is a unique point $v\in L^+$ on either side of the plane through the origin
containing
$u_0,u_1$ satisfying $\lambda (h(v),h(u_i))=\lambda _i$, for $i=0,1$.  The point $v$ depends continuously on $u_0,u_1$ and on $\lambda _0,\lambda
_1$.}

\vskip .1in 

\noindent {\bf Lemma~3} [7;Proposition 2.8] {\it 
Suppose that $u_0,u_1,u_2\in L^+$ are linearly independent, let $\gamma (u_i,u_j)$ denote the 
geodesic in ${\bf H}$ with ideal vertices given by the centers of $h(u_i)$ and $h(u_j)$, for $i\neq j$, and define
$$-\lambda _i^2= <u_j,u_k>,~~\alpha _i={{\lambda _i}\over{\lambda _j\lambda _k}},~~for~~\{ i,j,k\} =\{ 0,1,2\} .$$
Then $2\alpha _i$ is the hyperbolic length along the horocycle $h(u_i)$ between $\gamma (u_i,u_j)$ and $\gamma (u_i,u_k)$,
for $\{ i,j,k\} =\{ 0,1,2\}$.
}

\vskip .2in

\noindent {\bf Remark 1}~Consider an ideal quadrilateral $Q$ in ${\bf D}$ decorated so as to give four points in $L^+$.  The edges of $Q$ have well
defined lambda lengths, say $a,b,c,d$ in correct cyclic (clockwise) order about the boundary of $Q$.  Choose a diagonal of
$Q$, where the diagonal has lambda length $e$ and separates edges with lambda lengths $a,b$ from edges with lambda lengths $c,d$. The
other diagonal of $Q$ has its lambda length $f$ given by $$ef=ac+bd,$$
and we say that $f$ arises from $e$ by a {\it Ptolemy transformation} on the lambda lengths.  To see this, note that the formula for a Ptolemy
transformation is independent under scaling any of the four points in $L^+$, so we may alter the decoration and assume that the four points lie
in a common horizontal plane.  In this plane, the Minkoswki inner product induces a multiple of the usual Euclidean metric, and the intersection
of $L^+$ with this plane is a round circle.  The formula for the Ptolemy transformation thus follows from Ptolemy's classical formula on
Euclidean lengths of quadrilaterals that inscribe in a circle.

\vskip .2in

\noindent {\bf Remark 2}~Consider a decorated triangle, say with lambda lengths $x,y,z$ in the cyclic order about the boundary of the
triangle determined by an orientation, and define a 2-form $$\omega (x,y,z)=d{\ln}~x\wedge d{\ln}~y+d{\ln}~y\wedge d{\ln}~z+d{\ln}~z\wedge
d{\ln}~x.$$ A calculation shows that $\omega
(a,b,e)+\omega (c,d,e)=
\omega (b,c,f)+\omega (d,a,f)$, thus assigning a well-defined Ptolemy-invariant
 2-form to an oriented decorated ideal quadrilateral.

{\centerline{\epsffile{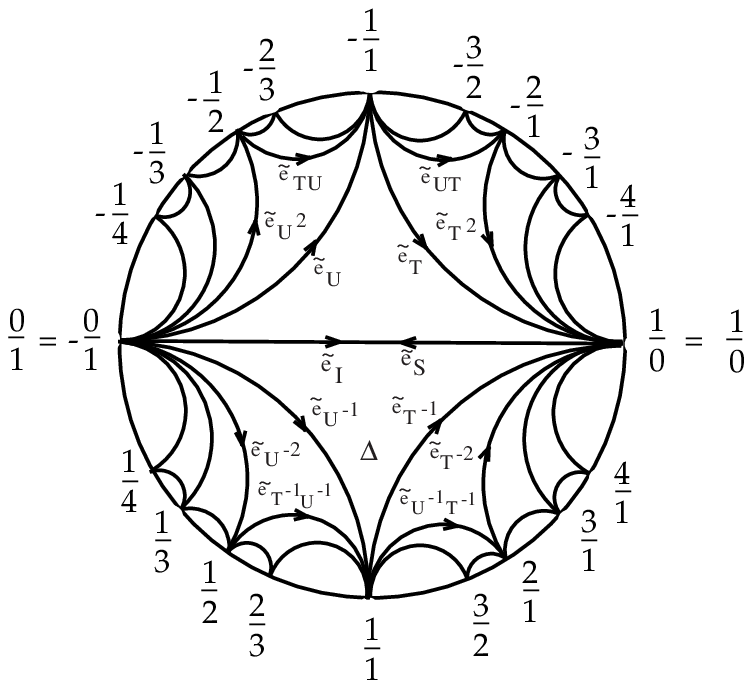}}}

\centerline{\bf Figure 2}

\vskip .2in

\noindent Regard the Poincar\'e disk as the open unit disk 
${\bf D}$ in the complex plane in the usual
way so that the unit circle $S^1$ is identified with the circle
at infinity, and let $\Delta$ denote
the ideal hyperbolic triangle with vertices $+1,-1, -\sqrt{-1}\in S^1$
as in Figure~2.
Let $\Gamma$ denote the group generated by reflections in
the sides of $\Delta$, and define
the {\it Farey tesselation} $\tau _*$ to be the full $\Gamma$-orbit
of the frontier of $\Delta$.  We refer to geodesics in $\tau _*$
as {\it edges} of $\tau _*$ and think of $\tau _*$ itself as a set of edges.
The ideal vertices of the edges of $\tau _*$ are naturally
identified with the set ${\bf Q}$
of all rational numbers including infinity, where for instance 
$+1,-1, -\sqrt{-1}\in S^1$
correspond respectively to $\infty={1\over 0}$, $0={0\over 1}$, $1={1\over 1}$
as illustrated in Figure~2.  Let ${\bf Q}
\subseteq S^1$
denote the corresponding countable dense subset of $S^1$
which we refer to simply as the set of {\it rational points} of $S^1$. 
Define the {\it distinguished oriented
edge} or {\it doe} of the Farey tesselation to be the
oriented edge from ${0\over 1}$ to ${1\over 0}$.

\vskip .1in

\noindent The {\it modular group} $PSL_2=PSL_2({\bf Z})$ of integral fractional
linear transformations is the subgroup of $\Gamma$
consisting of compositions of an even number of reflections, and
$PSL_2$ acts simply transitively on the set of orientations
on the edges of $\tau _*$.  The assignment 
$$
e_A = A(doe),~{\rm for}~A\in PSL_2,$$
establishes a bijection between $
PSL_2$ and the set of oriented
edges of $\tau _*$ as illustrated in Figure~2.  In particular,
the doe of $\tau _*$ is $e_I$, where 
$I$ denotes the identity of $PSL_2$.

\vskip .1in  

\noindent We adopt the standard notation 
$$S=\pmatrix{0&-1\cr 1&~~0\cr},~~T=
\pmatrix{ 1&1\cr 0&1},~~U=\pmatrix {1&0\cr 1&1\cr}$$
for certain elements of $PSL_2$,
where $S$ is involutive and 
fixes the unoriented edge of $\tau _*$ underlying the doe
while changing its orientation, and $U$ (respectively  $T$) is 
the parabolic transformation
with fixed point ${0\over 1}$ (respectively ${1\over 0}$)
which cyclically permutes the incident edges of $\tau _*$ in
the counter-clockwise sense about ${0\over 1}$ (respectively
the clockwise sense about ${1\over 0}$).  
In fact, $U^{-1}=STS,~T^{-1}=SUS$, and any two
of $S,T,U$ generate $PSL_2$.  

\vskip .1in

\noindent We shall also require the full {\it  M\"obius group} 
$M\ddot ob= PSL_2({\bf R})\supseteq PSL_2({\bf Z})=PSL_2$
consisting of all real fractional linear transformations.

\vskip .3in

\noindent {\bf 2. Punctured surfaces}

\vskip .2in

\noindent Let $F=F_g^s$ denote a fixed smooth surface of genus $g$ with $s\geq 1$ punctures, where $2-2g-s<0$.

\vskip .1in

\noindent Choose any base-point to determine the fundamental group $G$ of $F$, and consider the space $Hom'(G,M\ddot ob)$
of all discrete and faithful representations $\rho :G\to M\ddot ob$ so that no holonomy $\rho (\gamma )$ is
elliptic for $\gamma\in G $, and the holonomies around the punctures of
$F$ are parabolic.  Define the {\it Teichm\"uller space}  $$T(F)~=~{\rm Hom}'(G,M\ddot ob)/M\ddot ob,$$ where $M\ddot ob$ acts by
conjugacy.  

\vskip .1in

\noindent If $\rho\in Hom '(G,M\ddot ob)$, then ${\bf D}/\rho (G)$ induces a complete finite-area hyperbolic structure on $F$, whose
punctures are in one-to-one correspondence with the $G$-orbits of the set of fixed points of parabolic elements of $\rho (G)$.

\vskip .1in

\noindent 
Define the {\it decorated Teichm\"uller space} $\tilde T(F)\to T(F)$ of $F$ to be the trivial $R^s_{>0}$-bundle, where the fiber over a point is the
set of all $s$-tuples of horocycles, one horocycle about each puncture of $F$ parametrized by hyperbolic length

\vskip .1in

\noindent By an {\it arc family} in $F$, we mean the isotopy class of a family of essential arcs disjointly
embedded in
$F$ connecting punctures, where no two arcs in a family may be homotopic rel punctures.  If $\alpha$ is a maximal arc family, so that each component
of
$F-\cup \alpha$ is a triangle, then we say that $\alpha$ is an {\it ideal triangulation} of $F$.

\vskip .1in

\noindent {\bf Theorem 4} [7;Theorem 3.1] {\it Fix an ideal triangulation $\tau$ of $F$.  Then 
the assignment of $\lambda$-lengths $\tilde{\cal T}(F)\to R^\tau_{>0}$ is a homeomorphism onto.}

\vskip .1in 

\noindent
{\bf Proof}~~We must describe an inverse to the mapping and thus give the construction of a decorated hyperbolic structure from 
an assignment of putative $\lambda$-lengths.  To this end, consider the topological universal cover $\tilde F$ of $F$ and the lift 
$\tilde\tau$ of $\tau$ to
$\tilde F$; to each component arc of $\tilde\tau$ is associated the lambda length of its projection. 

\vskip .1in   

\noindent The proof proceeds by
induction, and for the basis step, choose any triangle $\Delta _0'$ of
$\tilde\tau$ and any ideal triangle $\Delta _0$ in
${\bf H}$.  The ideal vertices of $\Delta_0$ determine three rays in $L^+$, so by Lemma~1, there are three well-defined points in $L^+$
realizing the putative $\lambda$-lengths on the edges of $\Delta _0'$.  (In effect, this basis step of normalizing a triangle ``kills'' the
conjugacy by the M\"obius group in the definition of Teichm\"uller space.)  Of course, the triple of points in $L^+$ corresponds by affine
duality to a triple of horocycles, one centered at each ideal vertex of $\Delta _0$, i.e., a ``decoration'' on $\Delta _0$.

\vskip .1in   

\noindent To begin the induction step, consider a triangle $\Delta _1'$ adjacent to $\Delta_0'$ across an arc in $\tilde\tau$.  The two ideal
points which $\Delta _0'$ and $\Delta _1'$ share have been lifted to $u,v\in L^+$ in the basis step, and we let $w\in L^+$ denote the lift of
the third ideal point of
$\Delta _0'$ and consider the plane through the origin determined by $u,v$.  According to Lemma~2, there is a unique lift $z\in L^+$ of the third
ideal point of
$\Delta_1'$ on the side of
this plane not containing the lift of $w$, where $z$ realizes the putative $\lambda$-lengths.  Again, $u,v,z$ gives rise via affine duality to
another decorated triangle $\Delta _1$ in ${\bf H}$ sharing one edge and two horocycles with $\Delta_0$.

\vskip .1in   

\noindent One continues in this manner serially applying Lemma~2 to produce a collection of decorated triangles pairwise sharing edges in ${\bf H}$,
where any two triangles have disjoint interiors (because of our choice of the side of the plane in Lemma~2).  Thus, the construction gives an
injection
$\tilde F\to {\bf H}$, and we next show that in fact this mapping is also a surjection.  To this end, note first that the inductive
construction has an image which is open in ${\bf H}$ by construction.  According to Lemma~3, there is some $\varepsilon >0$ so that each
horocyclic arc inside of each triangle has length at least $\varepsilon$; indeed, there are only finitely many values for such lengths
because the surface is comprised of finitely many triangles.  Thus, each application of the inductive step moves a definite amount along each
horocycle, and it follows that the construction has an image which is closed as well.  It follows from connectivity of ${\bf H}$ that $\tilde
F\to{\bf H}$ is surjective, and furthermore, one can see that $\tilde\tau$ is a tesselation of ${\bf H}$, i.e., a locally finite
collection of geodesics decomposing ${\bf H}$ into ideal triangles.

\vskip .1in   

\noindent Following Poincar\'e, the hyperbolic symmetry group of this tesselation is the required (normalized) Fuchsian group $G$ giving a point
of Teichm\"uller space, and the construction likewise provides a decoration on the quotient ${\bf H}/G$ as required.
~~~~\hfill{\it q.e.d.}

\vskip .2in

\noindent {\bf Remark 3}~One thinks of the choice of ideal triangulation as a choice of ``basis'' for the lambda length coordinates.  Formulas
for the change of basis are given by Ptolemy transformations, and this leads [7],[8] to a faithful representation of the mapping class group of $F$
as well as its universal generalization to ${\cal T}ess$. 
Furthermore as in Remark 2, the 2-form
$\omega=\sum
\omega (a,b,c)$ is invariant under this action and descends to the Weil-Petersson form on Riemann's moduli space, where the sum is over all
triangles complementary to any fixed ideal triangulation and the edges of the triangle have lambda lengths $a,b,c$ in correct cyclic order
determined by an orientation of $F$.  These two ingredients lead to a natural quantization of Teichm\"uller space, [3] in the Poisson quantization
and [4] in the symplectic quantization.

\vskip .3in

\noindent {\bf 3. Coordinates for circle homeomorphisms}\vskip .2in

\noindent Define a {\it tesselation} $\tau$ of the Poincar\'e disk ${\bf D}$ to be a countable locally
finite collection of hyperbolic geodesics in ${\bf D}$ each of whose complementary regions is an ideal triangle.
A {\it distinguished oriented edge} or {\it doe} of $\tau$ is the specification of an orientation on
one of the geodesics in $\tau$.
Each geodesic in $\tau$ has a pair of asymptotes in $S^1$, and we let $\tau ^0\subseteq S^1$ denote the collection of all 
such asymptotes of geodesics in $\tau$ and $\tau ^2$ denote the collection of all triangles complementary to $\cup\tau$.

\vskip .1in

\noindent Tesselations with doe are ``combinatorially rigid'' in the following sense.  Suppose that $\tau _1,\tau _2$ are each a tesselation
with doe, say the initial and terminal points of the doe in $\tau _i$ are $x_i\in\tau _i^0$ and $y_i\in\tau _i^0$, respectively, for
$i=1,2$.  There is a unique bijection $f:\tau _1^0\to \tau _2^0$ so that $f(x_1)=x_2,f(y_1)=y_2$, and whenever $x,y,z$ in correct cyclic order span
an oriented triangle in
$\tau _1^2$, then $f(x),f(y),f(z)$ in correct cyclic order also span an oriented triangle in $\tau _2^2$.  This mapping $f:\tau _1^0\to\tau _2^0$
is called the {\it characteristic mapping} of the pair of tesselations with doe.  In particular, we may fix $\tau _1=\tau _*$ to be the Farey
tesselation with doe defined in $\S$1, so $\tau _*^0={\bf Q}$.  We may thus define the characteristic mapping
$f_{\tau}:{\bf Q}\to\tau ^0$ of the tesselation $\tau=\tau _2$ with doe. 

\vskip .1in

\noindent  Define the set 
$${\cal T}ess'=\{{\rm tesselations~with~doe~of}~{\bf D}\} .$$
To define a topology on ${\cal T}ess'$, if $\tau$ is a tesselation with doe, then we may extend the
range of the characteristic mapping $f_\tau:{\bf Q}\to\tau^0\subseteq S^1$ to $S^1$.  The assignment $\tau\mapsto f_\tau$ determines an embedding
of ${\cal T}ess'$ into the function space $(S^1)^{{\bf Q}}$ with the compact-open topology (where ${\bf Q}$ is given its discrete topology),
and we endow ${\cal T}ess'$ with the subspace topology.

\vskip .1in

\noindent Define the topological group $Homeo_+=Homeo_+(S^1)$ to be the group of all orientation-preserving homeomorphisms of the circle
taken with the compact-open topology.  If $f\in Homeo_+$ and $e$ is any geodesic in ${\bf D}$, say with ideal points $x,y\in S^1$, then define
$f(e)$ to be the geodesic in ${\bf D}$ with ideal points $f(x),f(y)\in S^1$.  It is not difficult to see that if $\tau$ is a tesselation and
$f\in Homeo_+$, then $f(\tau )=\{ f(e):e\in\tau\}$ is also a tesselation.  Since a doe on $\tau$ determines a doe on $f(\tau )$ in the natural
way, there is thus an action of $Homeo_+$ on ${\cal T}ess'$.

\vskip .2in

\noindent {\bf Theorem 5}~~[8;Theorem~2.3]~~\it The mapping
$$\eqalign{
Homeo_+&\to{\cal T}ess'\cr
f&\mapsto f(\tau _*)\cr
}$$
is a homeomorphism onto.\rm

\vskip .2in

\noindent {\bf Proof}~Injectivity follows from the fact that a homeomorphism is uniquely determined by its values on a dense set.
For surjectivity, consider any tesselation with doe
$\tau$.  Using the fact that
${\bf Q}$ and
$\tau ^0$ are dense in
$S^1$ and
the characteristic mapping $f_\tau$ is order-preserving by construction, a standard point-set topology argument show that there is a unique
orientation-preserving homeomorphism $f_\tau:S^1\to S^1$ which restricts to the characteristic mapping.  Both mappings $f\mapsto f(\tau )$ and
$\tau\mapsto f_\tau$ are continuous by construction.~~~~~\hfill{\it q.e.d.}

\vskip .2in

\noindent There is the natural diagonal left action of the group $M\ddot ob$ on $(S^1)^{{\bf Q}}$, which induces a left action of $M\ddot ob$ on
the subspace ${\cal T}ess'$, and we finally define the {\it universal Teichm\"uller space}
$${\cal T}ess=M\ddot ob \backslash {\cal T}ess'\approx M\ddot ob\backslash Homeo_+$$
to be the orbit space with the quotient topology.

\vskip .1in

\noindent A {\it decoration} on a tesselation $\tau$ is the specification of
horocycles in ${\bf D}$, one hororcyle centered at each point of $\tau ^0$.
Via the affine duality discussed in $\S$1, the characteristic mapping $f_\tau:{\bf Q}\to S^1$ on a decorated tesselation $\tau$ with doe
extends to a mapping $g_\tau :{\bf Q}\to L^+$.  The image $g_\tau({\bf Q})$ is automatically ``radially dense'' in $L^+$ in the sense that
$\Pi (g_\tau({\bf Q}))$ is a dense subset of $S^1$, where $\Pi :L^+\to S^1$ is the natural projection.
Define 
$$\widetilde{{\cal T}ess}'=\{ {\rm decorated~tesselations}~\tau ~{\rm with~doe}: g_\tau ({\bf Q})~{\rm is~discrete~in}~L^+\} .$$
The Hausdorff topology on the set of all closed subsets of $L^+$ induces a subspace topology on the set of all discrete subsets of $L^+$, and this
in turn induces a compact-open topology on $\widetilde{{\cal T}ess}'$.
There is again a diagonal left action of $M\ddot ob$ by Minkowski isometries on $\widetilde{{\cal T}ess}'$, and the {\it decorated universal
Teichm\"uller space} is finally defined to be the topological quotient
$$\widetilde{{\cal T}ess}=M\ddot ob\backslash\widetilde{{\cal T}ess}'.$$
There is the natural forgetful map $\widetilde{{\cal T}ess}\to{\cal T}ess$, which is evidently continuous.

\vskip .1in

\noindent Given a decorated tesselation $\tilde\tau$ with doe and $e\in\tau _*$, there is the corresponding lambda length
of the decorated geodesic in $\tilde\tau$ with underlying geodesic $f_\tau (e)$, where $f_\tau$ is the characteristic mapping of $\tau$.  Thus,
lambda lengths naturally determine an element of ${\bf R}_{>0}^{\tau _*}$.  

\vskip .2in

\noindent{\bf Theorem 6}~~[8:Theorem~3.1]~~\it The assignment of lambda lengths determines an embedding
$$\widetilde{{\cal T}ess}\to {\bf R}_{>0}^{\tau _*}$$
onto an open set, where ${\bf R}_{>0}^{\tau _*}$ is given the weak topology (compact-open on ${\bf R}_{>0}^{\tau _*}$ with $\tau _*$ discrete).
Thus, $\widetilde{{\cal T}ess}$ inherits the structure of a Fr\'echet manifold.\rm

\vskip .2in

\noindent{\bf Proof}~
We say that an element $\tau \in{\cal T}ess'$ is {\it normalized} provided that $\{\pm 1\}\subseteq \tau ^0$, the doe of $\tau $ runs
from $-1$ to $+1$, and the triangle in $\tau^2$ lying to the right of the doe coincides with the triangle spanned by $-1,+1,-\sqrt{-1}\in S^1$.
Since $M\ddot ob$ acts three-effectively on $S^1$ and the value of a M\"obius transformation at three points of $S^1$ determines it uniquely, each
$M\ddot ob$-orbit on ${\cal T}ess'$ admits a unique normalized representative.  ${\cal T}ess$ is thus canonically identified with the collection of
all normalized tesselations.  (Again, we have ``killed'' the M\"obius group by normalization.) 

\vskip .1in

\noindent To define a left inverse to the assignment $\lambda\in{\bf R}_{>0}^{\tau _*}$ of lambda lengths, use Lemma~1 to uniquely lift the vertices
$\pm 1,-\sqrt{-1}$ of the triangle of $\tau _*$ to the right of the doe to points in the rays in $L^+$ lying over these vertices realizing
the lambda lengths.  As in the proof of Theorem~4, we may then uniquely extend using Lemma~2 to a function $g:{\bf Q}\to L^+$ realizing the
lambda lengths.  

\vskip .1in

\noindent If $g({\bf Q})\subseteq  L^+$ is radially dense, then the order-preserving mapping ${\bf Q}\to L^+\to S^1$ interpolates a
unique homeomorphism $f:S^1\to S^1$ as before.  One can always produce a discrete decoration, say by taking the point $f(p)\in L^+$ to have
height $i$ in ${\bf R}^3$ if $p\in{\bf Q}$ is of Farey generation $i$. 

\vskip .1in

\noindent It follows that the mapping $\widetilde{{\cal T}ess}\to {\bf R}_{>0}^{\tau _*}$ is indeed injective.  Continuity follows from the
definition of the topologies, and openness of the image follows from the construction.~~~~~\hfill{\it q.e.d.}

\vskip .2in

\noindent Recall [1] that a {\it quasi-symmetric} homeomorphism of $S^1$ is the restriction to $S^1$ of a quasi-conformal homeomorphism
of
${\bf D}$, and let $Homeo_{qs}\subseteq Homeo_+$ denote the subspace of all quasi-symmetric homeomorphisms of the circle.  {\it Bers' universal
Teichm\"uller space} [2] is the quotient $$M\ddot ob\backslash Homeo_{qs}\subseteq M\ddot ob\backslash Homeo_+\approx{{\cal T}ess}$$ and is highly
studied.  

\vskip .1in

\noindent As is usual in these circumstances, it is difficult to explicitly characterize the image $\widetilde{{\cal T}ess}\subseteq{\bf
R}_{>0}^{\tau _*}$. On the other hand, there are the following useful characterizations of subsets of $\widetilde{{\cal T}ess}\subseteq
{\bf R}_{>0}^{\tau
_*}$.

\vskip .1in

\noindent We say that $\lambda\in {\bf R}_+^{\tau
_*}$ is {\it pinched} provided there is some real number $K>1$ so that
$$K^{-1}\leq \lambda (e)\leq K,~{\rm for~each}~e\in\tau _*.$$

\vskip .2in

\noindent {\bf Theorem~7}~[8;Theorem 6.3] \it ~\it If $\lambda \in{\bf R}^{\tau _*}_{>0}$ is pinched,
then there is a decorated tesselation whose lambda
lengths are given by $\lambda$, and the corresponding subset of $L^+$ is discrete and radially dense.\rm

\vskip .2in

\noindent {\bf Theorem 8}~[8;Theorem 6.4 (joint with Sullivan)]
~\it If $\lambda \in{\bf R}_{>0}^{\tau _*}$ is pinched,
then the corresponding homeomorphism
of the circle is quasi-symmetric.\rm
\vskip .1in

\leftskip=0ex\rightskip=0ex\rm

\noindent In particular, consider any decoration on a marked punctured Riemann surface $F$ uniformized by $G<M\ddot ob$.  Choose an ideal
triangulation of
${\bf D}/G$, and lift it to a tesselation $\tau$ of ${\bf D}$ which inherits a $G$-invariant decoration.  Choose a doe of $\tau$ to determine a
point of
$\widetilde{{\cal T}ess}$.  The lambda lengths in $F$ lift to $G$-invariant lambda lengths on $\tau$, and they are pinched since they take only
finitely many values.  Furthermore, if $G < PSL_2$ is finite-index and free of elliptics and $\phi:{\bf D}\to{\bf D}$ is any
$G$-invariant quasi-conformal map conjugating $G$ to an isomorphic group, then we claim that the boundary values of $\phi$ satisfy the smoothness
conditions above.  To see this, conjugate in domain and range so that corresponding parabolic covering transformations are each given in upper
halfspace by
$z\mapsto z+1$ (thus, destroying the normalization in $\widetilde{{\cal T}ess}$), so the monotone function $\phi (t)$ nearly agrees with the integral
part of
$t$.  It follows directly that
$\phi (t)$ is differentiable at each point of ${\bf Q}$, and the derivatives at points of ${\bf Q}$
are uniformly near unity.  (Compare with [13].)

\vskip .1in

\noindent In contrast, a quasi-symmetric map $\phi :S^1\to S^1$ arising from pinched lambda lengths need not have these differentiability properties
at
${\bf Q}$.  To see this, use the formula in the Introduction for lambda lengths in the upper halfspace model to produce pinched
lambda lengths so that the two one-sided derivatives at infinity disagree
or even fail to exist.

\vskip .1in

\noindent This corrects the second part of Theorem~6.4 from [8].  For the corresponding result in the setting of the solenoid, see Proposition~12 in $\S$4.

\vskip .4in

\noindent{\bf 4. Coordinates for the solenoid}

\vskip .2in

\noindent Let $G<PSL_2$ be any finite-index subgroup, and choose a base-point in the quotient surface or orbifold
$M={\bf D}/G$; in particular, for $G=PSL_2$, $M$ is the orbifold modular curve.  Consider the category ${\cal C}_M$ of all finite-sheeted orbifold
covers
$\pi: F\to M$, where $F$ is a punctured Riemann surface.  ${\cal C}_M$ is a directed set, where $\pi _1\leq \pi _2$ if there is a finite-sheeted
unbranched cover $\pi _{2,1}:F_2\to F_1$ of Riemann surfaces so that the following diagram commutes:

\vskip .1in

\settabs 10\columns

\+&&$F_2$&${{\pi_{2,1}}\atop\rightarrow}$&$F_1$\cr
\+\cr
\+&&~~$\pi _2\searrow$&~~~~$\swarrow\pi _1$\cr
\+\cr
\+&&&M\cr
\vskip .1in

\noindent In other words by covering space theory, if $\Gamma _i<G<PSL_2$ uniformizes $F_i$ for $i=1,2$, then $\pi _1\leq \pi _2$ if and only if
$\Gamma _1 $ is a finite-index subgroup of $\Gamma _2$. 

\vskip .1in

\noindent A topological space, the {\it punctured solenoid}, is defined in analogy to [12] to be the inverse limit 
$${\cal H}_M=\lim _\leftarrow{\cal C}_M;$$
a point of ${\cal H}_M$ is thus determined by choices of points $y_i\in F_i$ for each cover $\pi _i:F_i\to M$, where
the choices are ``compatible'' in the sense that  if $\pi _1\leq \pi _2$, then we have in the notation above $\pi _{2,1}(y_2)=y_1$.

\vskip .1in

\noindent Since punctured surface groups are cofinal in the set of punctured orbifold groups, we could have equivalently considered the category of
orbifold covers of $M$ in the definition of ${\cal H}_M$.  Furthermore, if $\Gamma < G$ is of finite-index, then ${\cal H}_\Gamma$ is naturally
homeomorphic to ${\cal H}_G$, and we may thus think of {\sl the} punctured solenoid ${\cal H}={\cal H}_{PSL_2}$

\vskip .1in

\noindent One can from first principals develop the Teichm\"uller theory of ${\cal H}$ along classical lines [12], [5], [11].  Instead, 
we next introduce an explicit space homeomorphic to
${\cal H}$ following [5], and we shall then {\sl define} the Teichm\"uller space representation theoretically in analogy to punctured surfaces
in $\S$1.

\vskip .1in

\noindent $G$ has characteristic subgroups  
$$G_N=\cap\{ \Gamma < G:[\Gamma :G]\leq N\},$$
for each $N\geq 1$, and these are nested $G_N < G_{N+1}$.  Define a metric 
$$\eqalign{
G\times G&\to{\bf R},\cr
\gamma\times\delta&\mapsto {\min}\{ {1\over N}:\gamma\delta ^{-1}\in G_N\} ,\cr
}$$
and define the {\it pro-finite completion} $\hat G$ of $G$ as a space to be the metric completion
of $G$, i.e., suitable equivalence classes of Cauchy sequences in $G$.  Termwise multiplication
of Cauchy sequences gives $\hat G$ the structure of a topological group, and termwise multiplication by $G$ gives a continuous
action of $G$ on $\hat G$.

\vskip .1in

\noindent For any sub-group $G<PSL_2$ of finite-index, we may define the quotient
$${\cal H}_G={\bf D}\times_G\hat G=({\bf D}\times \hat{G})/G,$$
where $\gamma \in G$ acts by
$$\eqalign{
\gamma :{\bf D}\times\hat{G}&\to{\bf D}\times\hat{G}\cr
(z,t)&\mapsto (\gamma z,t\gamma ^{-1}).\cr
}$$

\vskip .2in

\noindent {\bf Lemma 9}~[5]~~\it ${\cal H}$ is homeomorphic to ${\cal H}_G$ for any $G<PSL_2$ of
finite-index.  \rm

\vskip .2in

\noindent In particular, each path component, or ``leaf'', of ${\cal H}$ is homeomorphic to a disk (by residual finiteness of $G$), and each
leaf is dense in ${\cal H}$ (since $G$ is dense in $\hat{G}$).

\vskip .1in

\noindent Let us for definiteness simply fix $G=PSL_2$ and consider the collection $Hom'(G\times\hat{G},M\ddot ob)$ of all functions
$\rho :G\times\hat{G}\to M\ddot ob$ satisfying the following three properties:

\vskip .2in

\leftskip .3in

\noindent {Property 1:}~$\rho$ is continuous;

\vskip .2in

\noindent {Property 2:}~for each $\gamma _1,\gamma _2\in G$ and $t\in \hat {G}$, we have $$\rho(\gamma _1\circ\gamma _2,t)=\rho(\gamma _1,t\gamma
_2^{-1})\circ \rho(\gamma _2,t);$$

\vskip .2in

\noindent {Property 3:}~for every $t\in\hat{G}$, there is a quasi-conformal mapping $\phi _t:{\bf D}\to{\bf D}$ so that for every $\gamma\in G$, 
the following diagram commutes:

\vskip .2in

\vskip .2in

\settabs 7\columns
\+&$(z,t)$&$\mapsto$&$(\gamma z,t\gamma ^{-1})$\cr
\+\cr
\+&${\bf D}\times\hat{G}$&$\to$&${\bf D}\times\hat{G}$\cr
\+\cr
\+&$\phi _t\times{\rm id}$~~$\downarrow$&&$\downarrow ~~\phi_{t\gamma ^{-1}}\times{\rm id}$\cr
\+\cr
\+&${\bf D}\times\hat{G}$&$\to$&${\bf D}\times\hat{G}$\cr
\+\cr
\+&$(z,t)$&$\mapsto$&$(\rho(\gamma ,t)z,t\gamma ^{-1})$\cr

\vskip .2in

\noindent Furthermore, $\phi _t$ varies continuously in $t\in\hat{G}$ for the
common refinement of the $C^{\infty}$ topology of uniform convergence on compacta in ${\bf D}$ and the usual topology on
Bers' universal Teichm\"uller space $M\ddot ob \backslash Homeo_{qs}$ of the 
extension of $\phi _t$ to the circle at infinity.

\leftskip=0ex

\vskip .2in

\noindent  As to Property 1, notice that since $G$ is discrete, $\rho$ is continuous if and only if it is so in its second variable only.
Property 2 is a kind of homomorphism property of $\rho$ mixing the leaves; notice in particular that taking $\gamma _2=~I$ shows that
$\rho(I,t)=I$ for all $t\in\hat{G}$.  Property 3
mandates that for each
$t\in
\hat{G}$,
$\phi _t$ conjugates the standard action of $G$ on ${\bf D}\times\hat{G}$ at the top of the diagram to the action
$$\gamma _\rho:(z,t)\mapsto (\rho (\gamma ,t)z,t\gamma ^{-1})$$ at the bottom, and we let $G_\rho=\{ \gamma _\rho:\gamma\in G\}\approx G$.  
Notice that the action of $G_\rho$ on ${\bf D}\times\hat{G}$ extends continuously to an action on $S^1\times\hat{G}$.  We
finally define the solenoid (with marked hyperbolic structure)
$${\cal H}_\rho =({\bf D}\times _\rho\hat{G})=({\bf D}\times \hat{G})/G_\rho.$$
The collection $\phi _t$, for $t\in\hat G$, thus induces a homeomorphism ${\cal H}\to {\cal H}_\rho$. 

\vskip .1in

\noindent Define the group $Cont(\hat{G},M\ddot ob)$ to be the collection of all continuous maps $\alpha:\hat{G}\to M\ddot ob$, where the product
of two $\alpha ,\beta\in Cont(\hat{G},M\ddot ob)$ is taken pointwise $(\alpha\beta) (t)=\alpha (t)\circ\beta (t)$ in $M\ddot
ob$.
$\alpha\in Cont(\hat{G},M\ddot ob)$ acts continuously on $\rho\in Hom'(G\times \hat{G},M\ddot ob)$ according to
$$(\alpha\rho)(\gamma, t) = \alpha^{-1}(t\gamma ^{-1})\circ\rho(\gamma ,t)\circ \alpha (t).$$  

\vskip .2in

\noindent {\bf Theorem 10}~~[10]~~\it There is a natural homeomorphism of the Teichm\"uller space of the solenoid ${\cal H}$ with
$$T({\cal H})=Hom'(G\times\hat{G},M\ddot ob)/Cont(\hat{G},M\ddot ob).$$\rm

\vskip .2in

\noindent Rather than describe the proof here, we shall for simplicity simply take this quotient as the definition of the Teichm\"uller space
$T({\cal H})$.  Again with an eye towards simplicity here, rather than define punctures of solenoids intrinsically (as suitable equivalence
classes of ends of escaping rays), we can more simply proceed as follows.  Each $\phi _t:{\bf D}\to{\bf D}$
extends continuously to a quasi-symmetric mapping $\phi _t:S^1\to S^1$.  We say that a point $(p,t)\in S^1\times\hat{G}$ is a {\it
$\rho$-puncture} if
$\phi _t^{-1}(p)\in {\bf Q}$, and a {\it puncture} of ${\cal H}_\rho$ itself is a $G_\rho$-orbit of $\rho$-punctures.  A {\it $\rho$-horocycle} at a
$\rho$-puncture $(p,t)$ is the image under $\phi _t$ of a horocycle in ${\bf D}$ centered at $\phi _t^{-1}(p)$. 

\vskip .1in

\noindent  
A {\it decoration} on ${\cal H}_\rho$, or a ``decorated hyperbolic structure'' on ${\cal H}$, is a function $$\tilde\rho:G\times\hat{G}\times{\bf
Q}\to M\ddot ob\times L^+ ,$$ 
where $$\tilde\rho(\gamma ,t,q)=\rho(\gamma ,t)\times h(t,q)$$ with $\rho(\gamma ,t)\in Hom'(G\times\hat{G},M\ddot ob)$,
which satisfies the following conditions:

\vskip .2in

\leftskip .3in

\noindent {Property 4:}~for each $t\in\hat{G}$, the image $h(t,{\bf Q})\subseteq L^+$ is discrete and radially dense;

\vskip .1in

\noindent {Property 5:}~for each $q\in{\bf Q}$, the restriction $h(\cdot,q):\hat{G}\to L^+$ is continuous;

\vskip .1in

\noindent {Property 6:}~$\tilde\rho$ is $G$-invariant in the sense that if $\delta\in G$, then
$$\delta\circ \tilde\rho(\gamma,t,q)=\tilde\rho(\delta\gamma,t\delta ^{-1},\delta q),$$
where $\delta$ acts diagonally $\delta :(\gamma,q)\mapsto (\delta\gamma,\delta q)$ on $M\ddot ob\times L^+$ with
$\delta q$ the natural action of $G=PSL_2$ on $L^+$.

\vskip .2in

\leftskip=0ex 

\noindent Finally, let $Hom'(G\times\hat{G}\times{\bf Q},M\ddot ob\times L^+)$ denote the space of all decorated hyperbolic structures on ${\cal H}$ satisfying the
properties above, and define the {\it decorated Teichm\"uller space}
as the quotient
$$\tilde T({\cal H})=Hom'(G\times\hat{G}\times{\bf Q},M\ddot ob\times L^+)/Cont(\hat{G},M\ddot ob),$$
where $\alpha :\hat{G}\to M\ddot ob$ acts on $\tilde\rho$ by
$$(\alpha \tilde\rho)(\gamma,t,q)=\bigl(\alpha^{-1}(t\gamma^{-1})\circ\rho(\gamma ,t)\circ\alpha(t)\bigr)\times \bigl(h(t,\alpha(t) q)\bigr).$$
Forgetting decoration induces a continuous surjection $\tilde {T}({\cal H})\to T({\cal H})$.

\vskip .1in

\noindent There is a natural mapping $\lambda:\tilde T({\cal H})\to({\bf R}_{>0}^{\tau _*})^{\hat G}$ which assigns to
a function $\tilde\rho:G\times\hat{G}\times{\bf Q}\to M\ddot ob\times L^+$ the lambda length for the
$G_\rho$ metric of the $\rho$-horocycles determined by $h$ at the endpoints of the geodesic in ${\cal H}_\rho$ labeled by $\gamma$.

\vskip .2in

\noindent{\bf Theorem 11}~~[10] ~~\it The assignment of lambda lengths determines an embedding
$$\tilde{T}({\cal H})\to Cont(\hat{G}, {\bf R}_{>0}^{\tau _*})$$ onto an open set,
where we take the strong topology on ${\bf R}_{>0}^{\tau _*}$ and on $Cont(\hat{G},{\bf R}_{>0}^{\tau _*})$.\rm

\vskip .2in

\noindent{\bf Sketch of Proof}~~To prove the mapping is injective, we must again define the construction of decorated hyperbolic structure from 
a continuous $\lambda _t\in ({\bf R}_{>0}^{\tau _*})^{\hat G}$.    
 To this end, begin the definition of $\tilde\rho=\rho\times h$
on the triangle to the right of the doe in $\tau _*$ with lambda lengths given by $\lambda _t$.  As usual according to Lemma~1, we can uniquely
lift to a triple of points in $L^+$ lying over $\pm 1,-\sqrt{-1}$.  

\vskip .1in

\noindent 
It is easily seen from the identification of $G=PSL_2$ with the oriented edges
of $\tau _*$ that any $\gamma\in PSL_2$ can be written uniquely in the one of the following forms:

\vskip .2in

\leftskip .3in

\noindent ~~i)~$\gamma =I$;

\vskip .1in

\noindent ~ii)~$\gamma$ lies in the free semi-group generated by either $U,T$ or 
$U^{-1},T^{-1}$;

\vskip .1in

\noindent iii)~$\gamma$ arises from either i) or ii) by addition of the prefix $S$.

\vskip .2in

\leftskip=0ex

\noindent To define $\rho (\gamma ,t)\in M\ddot ob$, we shall specify an ideal triangle-doe pair $(\Delta ', e')$ in ${\bf D}$,
where $e''$ is oriented with $\Delta '$ to its right.
There is then
a unique
$\rho\in M\ddot ob$ mapping to the vertices of $\Delta '$ the vertices $\pm 1,-\sqrt{-1}$ of the triangle $\Delta$ to the right of the doe $e_I$ in
$\tau _*$ and mapping
$e_I$ to
$e'$.

\vskip .1in

\noindent Let us write $\gamma\in G$ in one of the forms i-iii) above.  Of course, if $\gamma =I$, then $\rho (\gamma ,t)=I$ as follows from the
functional equation in Property 2, and we take $(\Delta ',e')=(\Delta ,e_I)$.  

\vskip .1in

\noindent If $\gamma =S$, then let us begin with the lambda lengths $\lambda _t\in{\bf R}_{>0}^{\tau _*}$ on the edges of $\Delta$ and employ
Lemma~1 to uniquely realize a lift to $L^+$ of the vertices of this decoration on the triangle $\Delta$.  On the triangle to the left  of
the doe, consider the lambda lengths $\lambda _{t S^{-1}}=\lambda _{tS}$.  It need not be that $\lambda _t(e_I)=\lambda _{tS}(e_I)$,
and we re-scale, taking lambda lengths
$${{\lambda _t(e_I)}\over{\lambda _{tS}(e_I)}}~~\lambda _{tS}(\cdot )$$
on the edges $e_U,e_T$.  This defines lambda lengths on the edges of the quadrilateral in $\tau _*$ triangulated by the doe.  Again
using Lemma~2, there is a unique lift of $\sqrt{-1}$ to $L^+$ realizing the lambda lengths, and the projection of this point to $S^1$ is one of the
vertices of
$\Delta '$.  The other two vertices of $\Delta '$ are $\pm 1$, and the doe is $e'=e_S$, completing the definition in this case.  Notice
that this element of $M\ddot ob$ that maps $(\Delta ,e_I)\to(\Delta ',e_S)$ is necessarily involutive.

\vskip .1in

\noindent The case of any word in one of the semi-groups in ii) is handled by induction on the length, where for instance for any such $\gamma$ that
has a prefix
$U$, one begins from the lambda lengths $\lambda _t$ on $\Delta$ and uses the re-scaled lambda lengths
$\lambda _{tU^{-1}}$ on the edges $e_U,e_T$; the edge corresponding to $e_U$ is the doe of the first step.  

\vskip .1in

\noindent The remaining case iii) of
a word from one of the semi-groups with prefix $S$ is handled in exactly the same manner completing the definition of the
construction of $\rho :G\times\hat{G}\to M\ddot ob$.  The functional equation in Property 2 on $\rho$ follows by definition.  Furthermore,
since $\lambda _t\in {\bf R}_{>0}^{\tau _*}$ depends continuously on $t$ (because of the definition of the topology), $\rho (\gamma ,t)$ 
is also continuous in $t$ as required in Property 1; indeed, the entries of $\rho (\gamma ,t)\in M\ddot ob$ are algebraic function of finitely many
lambda lengths.

\vskip .1in

\noindent As for Property 3 in the definition of $T({\cal H})$, we claim that for each
$t\in\hat{G}$, $\lambda_t\in{\bf R}_{>0}^{\tau _*}$ is pinched if $\lambda\in Cont(\hat{G},{\bf R}_{>0}^{\tau _*})$.

\vskip .1in

\noindent To see this, note that the very definition of a continuous function $\lambda :\hat{G}\to{\bf R}_{>0}^{\tau _*}$ means that 
$\forall K ~\exists N ~\forall e\in\tau _*~\forall\gamma\in G_N$, we have
$$1+K^{-1}~~\leq ~~{{\lambda _t(e)}\over{\lambda _t(\gamma e)}}~~\leq 1+K.$$
Take say $K={1\over 2}$ and its corresponding $N$.  A fundamental domain for $G_N$ has only a finite collection of values of
lambda lengths, and any other lambda length is at most a factor of 3/2 times a lambda length this finite set, and at least a factor of 1/2
times a lambda length in this finite set.  $\lambda _t$ is therefore pinched, proving the claim.

\vskip .1in

\noindent It follows that for each $t\in\hat{G}$, $\lambda _t:\tau _*\to{\bf R}_{>0}$ is necessarily pinched.  By Theorem~8,
there is a corresponding quasi-conformal homeomorphism $\phi _t:{\bf D}\to{\bf D}$.  Commutativity of the diagram and continuity in Property 3
follow by construction, and this completes the proof that the function
$\rho$ constructed above lies in $Hom'(G\times\hat{G},M\ddot ob)$ and hence determines a point of $T({\cal H})$.

\vskip .1in

\noindent To define a decoration
on the $\rho$-punctures, each $\lambda _t$ determines a decoration on $\tau _*\times\{ t\}\subseteq{\bf D}\times\hat{G}$, as required. 
Property 4 is guaranteed by the claim and Theorem~7.  Property 5 holds as before in the strong sense that
the Euclidean coordinates of each $h(t,q)$ are algebraic functions of finitely many lambda lengths, and Property 6 holds by invariance of lambda
lengths under
$M\ddot ob$. ~~~~~\hfill{\it q.e.d.}

\vskip .2in

\noindent {\bf Proposition 12}~~[10]~~\it ~Suppose that $\lambda\in Cont(\hat{G},{\bf R}_{>0}^{\tau _*})$.  Then for each $t\in\hat {G}$,
$\lambda _t:\tau _*\to{\bf R}_{>0}$ corresponds to a quasi-conformal homeomorphism $\phi _t:{\bf D}\to {\bf D}$ whose quasi-symmetric extension
$\phi _t:S^1\to S^1$ is differentiable at each point of ${\bf Q}$ with derivative uniformly near unity.\rm

\vskip .2in

\noindent {\bf Proof}~As above, continuity of $\lambda _t$ in $t$ implies that each $\lambda _t$ is pinched, which gives a quasi-symmetric map
$\phi _t$, for
$t\in
\hat{G}$, by Theorem~8.  Using the upper halfspace model, normalize $\phi _t$ such that it fixes $0$ and $\infty$, whence the geodesics of $\tau _*$
that limit to
$\infty$ are mapped by $\phi _t$ onto geodesics which likewise limit to $\infty$.
Again by continuity, we conclude that for each $\epsilon >0$, $e\in\tau _*$ and $\gamma\in PSL_2$ fixing $\infty$, we have
$|\lambda _t(e)-\lambda _t(\gamma ^n e)|<\epsilon$ for $n$ sufficiently large.
The differences $a_n=\phi _t(n)-\phi _t(n-1)$ are then $\epsilon_1$ close using 
continuity of the assignment in Lemmas 1 and 2 of decorated ideal triangles given
lambda lengths.

\vskip .1in

\noindent We show that $\lim_{n\to\infty}{1\over n}~{\phi _t(n)}$ exists and is bounded, which proves the proposition.
To this end 
since $|a_i-a_{nk+i}|<\epsilon_1$ for all $i,k$, we find
$$
|{1\over n} (a_1+\cdots +a_n)-{1\over{nk}}(a_1+\cdots +a_{nk})|\leq
1/n~\sum_{i=1}^n |a_i-{1\over k}(a_i+a_{i+n}+\cdots +a_{i+n(k-1)})|
\leq \epsilon_1,
$$
and it follows that ${1\over n}~\phi _t(n)$ is a Cauchy sequence, as desired.
~~~~~\hfill{\it q.e.d.}

\vskip .2in

\noindent Differentibility at the rational points, which holds in the case of punctured surfaces (cf. the discussion following
Theorem~8) thus also holds for the solenoid according to Proposition~12.  Differentiability does not hold for general decorated tesselations
with pinched lambda lengths however (cf. the discussion and counter-example following Theorem~8).

\vskip .3in

\noindent {\bf 5. Concluding remarks}

\vskip .2in

\vskip .1in

\noindent In addition to the function $h:\hat G\times{\bf Q}\to L^+$ constructed in the proof of Theorem~11, there is another natural function
$h_1:\hat G\times{\bf Q}\to L^+$ defined as follows.  Begin with the ``standard'' decoration on
$\tau _*$ where all lambda lengths are
$\sqrt{2}$.  By Proposition~12, each $\phi _t$ is smooth with bounded derivative $\phi _t'(q)$ at each point of ${\bf Q}\subseteq S^1$.
In an upper halfspace model of ${\bf H}$ with the endpoint of the doe at infinity, scale the Euclidean diameter of the horocycle centered at
$q\in{\bf Q}$ by the absolute value of the derivative
$|\phi _t'(q)|$ to determine the diameter of the horocycle at $\phi _t(q)$.  This defines the function $h_1(t,q)$.  Notice that $h_1(t,{\bf Q})$
is again discrete and radially dense in $L^+$, but there seems to be no guarantee that the function $h_1$ discussed above satisfies Property 5.

\vskip .1in 

\noindent The representation-theoretic treatment of the Teichm\"uller theory of the punctured solenoid seems to us quite appealing.  For example,
the strong topology for the solenoid (in Theorem 11) in contrast to the weak topology for circle homeomorphisms (in Theorem~6) is interesting. 
Furthermore, one can naturally induce stronger or weaker  transverse structures in the $\hat{G}$ direction in ${\cal H}$ by imposing conditions
other than continuity on the lambda length functions in Theorem~11, and we wonder in particular what is the transverse regularity of the function $h_1$. 

\vskip .1in 

\noindent As was
mentioned in the Introduction, Theorem~11 is the first step of an ongoing program [10] to develop the decorated Teichm\"uller theory of the punctured
solenoid. Though there is an alternative to the construction of $\rho:G\times \hat{G}\to M\ddot ob$ in the proof of Theorem~11, the argument
given here bears a close relation to the treatment of broken hyperbolic structures in [6].

\vskip .1in

\noindent Lambda lengths enjoy the simple transformation property of Ptolemy transformations (cf. Remark~1), and furthermore, there is a simple
invariant two-form (cf. Remark~2).  These ingredients have been used [3],[4] to give a quantization of classical
Teichm\"uller theory (cf. Remark~3).  These same two ingredients persist for the universal Teichm\"uller theory as well as for
the punctured solenoid and might be used to quantize these Teichm\"uller theories as well.  

\vskip .1in

\noindent It is an open (but not centrally important) problem in number theory to calculate the index of $G_N$ in $G=PSL_2$, and an algorithm for its
calculation devolves to the ``cubic fatgraph enumeration problem'' for surfaces of fixed Euler characteristic arising as total spaces of 
covers of ${\bf D}/G$ of degree $N$.  More speculatively, there is a natural group [10] generated by the
$G_N$-equivariant Ptolemy moves for some $N$, which is closely related to the completions [9] of the universal Ptolemy group studied in the
context of Grothendieck absolute Galois theory.

\vfill\eject

\noindent {\bf Bibliography}

\vskip .2in

\noindent [1]~~L. V. Ahlfors, ``Lectures on Quasiconformal Mapping'', Wadsworth and Brooks Cole Advanced Books and Software, Monetery, CA (1987).

\vskip .1in

\noindent [2]~~L. Bers, ``Universal Teichm\"uller space'', in ``Analytic Methods in Mathematical Physics, Gordon and Breach, New York, (1968), 65-83.

\vskip .1in

\noindent [3]~~L. Chekhov and V. Fock, {``A quantum Techm\"uller space''}, {\it Theor. Math.
Phys.}, {\bf 120} (1999) 1245--1259; {\it Quantum mapping class group,
pentagon relation, and geodesics}
{\it Proc. Steklov Math Inst}, 149-163.

\vskip .1in

\noindent [4]~~R.~M.~Kashaev, {`` Quantization of Teichm\"uller spaces and the
quantum dilogarithm''}, {\it Lett. Math. Phys.}, {\bf 43}, No.~2,
(1998).

\vskip .1in

\noindent [5]~~C. Odden, ``Virtual automorphism group of the fundamental group of a closed surface'', Ph. D. Thesis, Duke University (1997).

\vskip .1in

\noindent [6]~~A. Papdopoulos and R. C. Penner, ``The Weil-Petersson K\"ahler form and affine foliations on surfaces'', to appear {\it Ann.
Glob. Anal. Geom.} (2004).

\vskip .1in

\noindent [7]~~ R. C. Penner, ``The decorated Teichm\"uller space of  punctured surfaces", 
{\it Comm. Math.  Phys.}  {\bf 113}   (1987),  299-339.

\vskip .1in

\noindent [8]~~---,``Universal constructions in Teichm\"uller theory",  
{\it Adv. Math.}  {\bf   98} (1993), 143-215.

\vskip .1in

\noindent [9]~~---,``The universal Ptolemy group and its completions'', Geometric 
Galois Actions II, London Math Society Lecture Notes {\bf 243}, 
Cambridge University Press (1997), eds. P. Lochak and L. Schneps.

\vskip .1in

\noindent [10]~--- and D. \v Sari\' c, ``Teichm\"uller theory of the punctured solenoid'', in preparation (2004).

\vskip .1in

\noindent [11]~~D. \v Sari\' c, ``On quasiconformal deformations of the universal hyperbolic solenoid'', preprint (2003).

\vskip .1in

\noindent [12]~~D. Sullivan, ``Linking the universalities of Milnor-Thurston, Feigenbaum and Ahlfors-Bers'', Milnor Festschrift, Topological methods
in modern mathematics (L. Goldberg and A. Phillips, eds.), Publish or Perish (1993), 543-563.

\vskip .1in

\noindent [13]~~P. Tukia, ``Differentiability and rigidity of M\"obius groups'', {\it Invent. Math.}~{\bf 82} (1985), 557-578.

\bye